\let\LocalVersion=Y
\let\MathTime=N
\let\MathPi=N

\if Y\LocalVersion
    \documentclass[10pt,twoside]{article}
\else 
    \documentclass{gtart}
\fi

\usepackage{amsmath}
\usepackage{amsthm}
\usepackage{amsfonts,amssymb}
\usepackage{graphicx}

\if Y\LocalVersion
    \if Y\MathTime
        \usepackage[nobold]{mathtime}
    \fi    
    \if Y\MathPi
        \usepackage[mathcal,mathbb]{mathpi}
    \fi    
    \usepackage{macros}
    \newcommand{\ack}{Acknowledgements.}
    
    \let\maketitlepage\maketitle
\else
    \newcommand{\ack}{Acknowledgements}
\fi    

\theoremstyle{plain}
\newtheorem{thm}{Theorem}
\newtheorem{cor}[thm]{Corollary}

\newtheorem{lem}[thm]{Lemma}
\newtheorem{prop}[thm]{Proposition}
\theoremstyle{definition}
\newtheorem{defn}[thm]{Definition}

\newcommand{\R}{\mathbb{R}}

\newcommand{\Q}{\mathbb{Q}}
\newcommand{\CP}{\mathbb{CP}}
\newcommand{\RP}{\mathbb{RP}}
\newcommand{\Z}{\mathbb{Z}}

\newcommand{\SO}{\mathit{SO}}
\newcommand{\so}{\mathfrak{so}}
\newcommand{\SU}{\mathit{SU}}
\newcommand{\A}{\mathbb{A}}

\makeatletter
\@ifpackageloaded{mathpi}{%
   
   }{}
\makeatother

\newcommand{\Hol}{\mathrm{Hol}}
\newcommand{\bonf}{\mathcal{B}}
\newcommand{\conf}{\mathcal{A}}
\newcommand{\G}{\mathcal{G}}
\newcommand{\Rep}{\mathcal{R}}
\newcommand{\g}{\mathfrak{g}}
\newcommand{\CS}{\mathrm{CS}}
\newcommand{\tr}{\mathrm{tr}}
\newcommand{\im}{\mathrm{im}}
\newcommand{\Aut}{\mathrm{Aut}}
\newcommand{\inc}{\iota}
\newcommand{\cE}{\mathcal{E}}
\newcommand{\grad}{\mathrm{grad}}
\newcommand{\sa}{a}
\renewcommand{\ss}{s}

\if Y\LocalVersion
    \title{Dehn surgery, the fundamental group and {\itshape SU}(2)}
    \author  {P. B. Kronheimer and T. S. Mrowka}
    \address  {Harvard University, Cambridge MA 02138\\
    Massachusetts Institute of Technology, Cambridge MA 02139}
\else
    \title{Dehn surgery, the fundamental group and {\itshape SU}(2)}
    \authors  {P. B. Kronheimer\\T. S. Mrowka}
    \address  {Harvard University, Cambridge MA 02138}
    \secondaddress  {Massachusetts Institute of Technology,
                     Cambridge MA 02139} 
    \email {kronheim@math.harvard.edu}
    \secondemail {mrowka@math.mit.edu}
\fi

\begin{abstract}
       Let $K$ be a non-trivial knot in the $3$-sphere and let $Y_{r}$ be the
        $3$-manifold obtained by surgery on $K$ with surgery-coefficient
        a rational number $r$. We show that there is a homomorphism
        from
        $\pi_{1}(Y_{r})$  to $SU(2)$ with
        non-cyclic image if $|r| \le 2$.
\end{abstract}

\begin{document}
\maketitlepage

\section{Introduction}

The main result of this paper, which is a companion to \cite{KM-PP},
is the following theorem.

\begin{thm}\label{thm:hRP3}
    Let $K$ be a non-trivial knot in $S^{3}$, and let $Y_{r}$ be the $3$-manifold
    obtained by Dehn surgery on $K$ with surgery-coefficient $r \in
    \Q$. If\/ $|r| \le 2$,  then
    $\pi_{1}(Y_{r})$ is not cyclic. In fact,
    there is a homomorphism $\rho : \pi_{1}(Y_{r}) \to
    \SU(2)$ with non-cyclic image. 
\end{thm}

The statement that $Y_{r}$ cannot have cyclic fundamental group was
previously known for all cases except $r=\pm 2$. The case $r=0$ is
due to Gabai \cite{Gabai3}, the case $r=\pm 1$
is the main result of \cite{KM-PP}, and the case that $K$ is a torus knot
is analysed for all $r$ in \cite{Moser}. All remaining cases follow from
the cyclic surgery theorem of Culler, Gordon, Luecke and Schalen
\cite{CGLS}.
It is 
proved in \cite{KMOS} that $Y_{2}$ cannot be homeomorphic to
$\RP^{3}$. If one knew that $\RP^{3}$ was the only closed
$3$-manifold with fundamental group $\Z/2\Z$ (a statement that is
contained in Thurston's geometrization conjecture), then the first
statement in the above theorem would be a consequence. The second
statement in the theorem
appears to sharpen the result slightly. In any event, we have:

\begin{cor}
    Dehn surgery on a non-trivial knot cannot yield a $3$-manifold with the
    same homotopy type as $\RP^{3}$. \qed
\end{cor}

The proof of Theorem~\ref{thm:hRP3} provides a verification of the Property P
conjecture that is independent of the results of the cyclic
surgery theorem of \cite{CGLS}.
Although the argument follows \cite{KM-PP} very closely, we shall avoid
making explicit use of instanton Floer homology and Floer's exact
triangle \cite{Floer2,Braam-Donaldson}. Instead, we rely on the
technique that
forms just the first
step of Floer's proof from \cite{Floer2}, namely the technique of
``holonomy perturbations'' for the instanton equations (see also the
remark following Proposition~16 in \cite{KM-PP}).

\subparagraph{\ack}
The first author was supported by NSF grant
DMS-0100771. The second author was supported by NSF grants 
DMS-0206485, DMS-0111298 and FRG-0244663. 

\section{Holonomy perturbations}

This section is a summary of material related to the ``holonomy
perturbations'' which Floer used in the proof of his surgery exact
triangle for instanton Floer homology \cite{Floer2}. Similar holonomy
perturbations were introduced for the $4$-dimensional
anti-self-duality equations in \cite{Donaldson-Orientations}; see also
\cite{Taubes-Cass}. Our exposition is
taken largely from \cite{Braam-Donaldson} with only small changes in
notation. Some of our gauge-theory notation is taken from
\cite{KM-structure}.

Let $Y$ be a compact, connected $3$-manifold, possibly with boundary.
Let $w$ be a unitary line bundle on
$Y$, and let $E$ be a unitary rank 2 bundle equipped with an
isomorphism
\[
         \psi : \det(E) \to w.
\]
Let $\g_{E}$ denote the bundle whose sections are the traceless, skew-hermitian
endomorphisms of $E$, and let $\conf$ be the affine space of $\SO(3)$
connections in $\g_{E}$.  Let $\G$ be the gauge group of unitary
automorphisms of $E$ of determinant $1$  (the automorphisms that
respect $\psi$). We write $\bonf^{w}(Y)$ for the quotient space
$\conf/\G$. A connection $A$, or its gauge-equivalence class $[A] \in
\bonf^{w}(Y)$, is \emph{irreducible} if the stabilizer of $A$ is the
group $\{\pm 1\} \subset \G$, and is otherwise \emph{reducible}. The
reducible connections are the ones that preserve a decomposition of
$\g_{E}$ as $\R\oplus L$, where $L$ is an orientable $2$-plane bundle;
these connections have stabilizer either $S^{1}$ or (in the case
of the product connection) the group $\SU(2)$.

\begin{defn}
    We write $\Rep^{w}(Y)\subset \bonf^{w}(Y)$ for the space of
    $\G$-orbits of flat connections:
    \[
              \Rep^{w}(Y) = \{ \,[A] \in \bonf^{w}(Y) \mid F_{A} = 0
              \,\}. 
    \]
    This is the \emph{representation variety} of flat connections with
    determinant $w$. \qed
\end{defn}

We have the following straightforward fact:

\begin{lem}
    The representation variety $\Rep^{w}(Y)$ is non-empty if and only
    if $\pi_{1}(Y)$ admits a homomorphism $\rho :
    \pi_{1}(Y) \to \SO(3)$ with $w_{2}(\rho) = c_{1}(w)$ mod $2$. The
    representation variety contains an irreducible element if and only
    if there is such a $\rho$ whose image is not cyclic.

    If $c_{1}(w)=0$ mod $2$, then $\Rep^{w}(Y)$ is isomorphic to the
    space of homomorphisms $\rho:\pi_{1}(Y)\to \SU(2)$ modulo the
    action of conjugation.\qed
\end{lem}

Suppose now that $Y$ is a closed oriented $3$-manifold.
The flat connections $A\in \conf$ are the critical points of the
Chern-Simons function
\[
\begin{gathered}
\CS : \conf \to \R ,\\
\CS(A) = \frac{1}{4}\int_{Y} \tr \bigl((A- A_{0}) \wedge (F_{A}+
F_{A_{0}})\bigr),
\end{gathered}
\]
where $A_{0}$ is a chosen reference point in $\conf$, and $\tr$
denotes the trace on $3$-by-$3$ matrices.  We define a class of
perturbations of the Chern-Simons functional, the \emph{holonomy}
perturbations.

Let $D$ be a
compact $2$-manifold with boundary, and let $\inc :
S^{1}\times D\hookrightarrow Y$. Choose a trivialization of $w$ over
the image of $\inc$. With this choice, each connection $A\in \conf$ gives
rise to a unique connection $\tilde A$ in $E|_{\im(\inc)}$ with the
property that $\det(\tilde A)$ is the product connection in the
trivialized bundle $w|_{\im(\inc)}$. Thus $\tilde{A}|_{\im(\inc)}$ is an
$\SU(2)$ connection. Given a smooth $2$-form $\mu$ with
compact support in the interior of $D$ and integral $1$, and given a smooth
class-function
\[
          \phi : \SU(2) \to \R,
\]
we can construct a function
\[
           \Phi : \conf \to \R
\]
that is invariant under $\G$ as follows.  For each $z \in D$, let
$\gamma_{z}$ be the loop $t \mapsto \inc(t,z)$ in $Y$, and let
$\Hol_{\gamma_{z}}(\tilde{A})$ denote the holonomy of $\tilde{A}$ along $\gamma_{z}$,
as an automorphism of the fiber  $E$ at the point $y = \inc(0,z)$.  The
class-function $\phi$ determines also a function on the group of
determinant-1 automorphisms of the fiber $E_{y}$, and we set
\[
               \Phi(A) = \int_{D} \phi(\Hol_{\gamma_{z}}(\tilde{A}))\mu(z).
\]

One can write down the equations for a critical point $A$ of the
function $\CS + \Phi$ on $\conf$. They take the form
\[
           F_{A} = \phi'(H_{A})\mu_{Y},
\]
where $H_{A}$ is the section of the bundle $\Aut(E)$ over $\im(\inc)$
obtained by taking holonomy around the circles, $\phi'$ is the
derivative of $\phi$, regarded as a map from $\Aut(E)$ to $\g_{E}$,
and $\mu_{Y}$ is the $2$-form on $Y$ obtained by pulling back $\mu$ to
$S^{1}\times D$ and then pushing forward along $\inc$. (See
\cite{Braam-Donaldson}.)

\begin{defn}
        Given $\inc$ and $\phi$ as above, we write
        \[
                  \Rep^{w}_{\inc,\phi}(Y) = \{ \, [A] 
                  \in \bonf^{w}(Y) \mid F_{A} = \phi'(H_{A})\mu_{Y}
              \,\}. 
        \]
        This is the \emph{perturbed representation variety}. \qed
\end{defn}

Now specialize to the case that $D$ is a disk, so $\inc$ is an
embedding of a solid torus. Let
\[
                C = Y \setminus \im(\inc)^{\circ}
\]
be the complementary manifold with torus boundary. Let $z_{0} \in
\partial D$ be a base-point, and let $a$ and $b$ be the oriented circles in
$\partial C$ described by
\begin{equation}\label{eq:ab-curves}
\begin{gathered}
        a =  \inc \bigl(S^{1} \times \{z_{0}\}\bigr) \\
        b    = \inc \bigl(\{ 0 \} \times \partial D\bigr).
\end{gathered}
\end{equation}
These are the ``longitude'' and ``meridian'' of the solid torus.
We continue to suppose that $w$ is trivialized on $\im(\inc)$ and
hence on $\partial C$. So the restriction of $E$ to $\partial C$ is
given the structure of an $\SU(2)$ bundle.
Given a connection $A$ on $\g_{E}$ that is flat
on $\partial C$, let $\tilde A$ be the corresponding flat $\SU(2)$
connection in $E|_{\partial C}$.  One can choose a determinant-1
isomorphism between the fiber of $E$ at the basepoint $\inc(0,z_{0})$
so that the holonomies of $\tilde A$ around $a$ and $b$ become
commuting elements of $\SU(2)$ given by
\[
\begin{aligned}
\Hol_{a}(\tilde A) &= 
        \begin{bmatrix}
                e^{i\alpha} & 0 \\
                0 & e^{-i\alpha}
        \end{bmatrix}
        \\
        \Hol_{b}(\tilde A) &= 
        \begin{bmatrix}
                e^{i\beta} & 0 \\
                0 & e^{-i\beta}
        \end{bmatrix}.
        \end{aligned}
\]
The pair $(\alpha(A),\beta(A)) \in \R^{2}$ is determined by $A$ up to the
ambiguities
\begin{enumerate}
        \item adding integer multiples of $2\pi$ to $\alpha$ or
        $\beta$;
        \item replacing $(\alpha,\beta)$ by $(-\alpha,-\beta)$.
\end{enumerate}

\begin{defn}
        Let $S \subset \R^{2}$ be a subset of the plane with the
        property that $S + 2\pi \Z^{2}$ is invariant under $s
        \mapsto -s$. Define the set \[ \Rep^{w}(C \mid S) \subset
        \Rep^{w}(C)
        \]
        as
        \[
             \Rep^{w}(C \mid S) =   \{ \,[A] \in \Rep^{w}(C) \mid
                   (\alpha(A),\beta(A)) \in S + 2\pi \, \Z^{2} \,\},
        \]
        where $(\alpha(A),\beta(A))$ are the longitudinal and meridional 
        holonomy parameters, determined up to the
        ambiguities above. \qed
\end{defn}

One should remember that the choice of trivialization of $w$ on
$\im(\inc)$ is used in this definition, and in general the set we have
defined will depend on this choice.

A class-function $\phi$ on $\SU(2)$ corresponds to a function $f:
\R\to\R$ via
\[
           f(t) = \phi\left(
           \begin{bmatrix}
              e^{it} & 0 \\
                0 & e^{-it}
            \end{bmatrix}
           \right).
\]
The function $f$ satisfies $f(t) = f(t + 2\pi)$ and $f(-t) = f(t)$.
The following observation of Floer's is proved as Lemma~5 in
\cite{Braam-Donaldson}.

\begin{lem}\label{lem:one-to-one-phi}
         Let $f:\R\to\R$ correspond to $\phi$ as above. Then
        restriction from $Y$ to $C$ gives rise to a bijection
        \[
               \Rep^{w}_{\inc,\phi}( Y ) \to
                       \Rep^{w}\bigl( \,C \mid \beta = -f'(\alpha)
                       \,\bigr).
        \]
        \qed
\end{lem}

We also have the straightforward fact:

\begin{lem}\label{lem:g}
        If $g : \R\to\R$ is a smooth odd function with period $2\pi$,
        then there is a class-function $\phi$ on $\SU(2)$ such that
        the corresponding function $f$ satisfies $f' = g$. \qed
\end{lem}

\section{Removing flat connections by perturbation}

Let us now take the case that $Y$ is a homology $S^{1}\times S^{2}$,
and let $w \to Y$ be a line-bundle with $c_{1}(w)$ a generator for
$H^{2}(Y;\Z) = \Z$.  Let $N \hookrightarrow Y$ be an embedded solid
torus whose core is a curve representing a generator of $H_{1}(Y;\Z)$,
and let $C$ be the manifold with torus boundary
\[
             C = Y \setminus N^{\circ}.
\]
By a ``slope'' we mean an isotopy class of essential closed curves on
the torus $\partial C$. For each slope $s$, let $Y_{s}$ denote the
manifold obtained from $C$ by Dehn filling with slope $s$: that is,
$Y_{s}$ is obtained from $C$ by attaching a solid torus in such a way
that curves in the class $s$ bound disks in the solid torus.

Parametrize $N$ by a map $\inc: S^{1} \times D^{2} \to N$. Let $a$ and
$b$ be the curves \eqref{eq:ab-curves} on $\partial N$. The Dehn
filling $Y_{b}$ on the slope represented by $b$ is just $Y$. The
manifold $Y_{a}$ has $H_{1}(Y_{a};\Z) = 0$. Let $\ss$ be the slope
\[
        \ss = [pa + qb],
\]
where $p$ and $q$ are coprime and both positive

\begin{prop}\label{prop:1/q-triangle}
        Let 
        $\ss$ be as above, and suppose
        \[
                  p/q \le 2.
        \]
        Suppose that neither
        $\pi_{1}(Y_{a})$ nor $\pi_{1}(Y_{\ss})$ admits a 
        homomorphism to $\SU(2)$ with non-cyclic image.
        Then there is a
        holonomy-perturbation $(\inc,\phi)$ for the manifold $Y$ such
        that the perturbed representation variety
        $\Rep^{w}_{\inc,\phi}(Y)$ is empty.
\end{prop}

\begin{proof}
Fix a trivialization $\tau$ of $w$ over $N$.
At this stage the choice is immaterial, because any two choices
differ by an automorphism of $w$ that extends over all of $Y$. Write
\[
\begin{aligned}
Y_{\sa} &= C \cup N_{\sa}, \\
Y_{\ss} &= C \cup N_{\ss}, \\
\end{aligned}
\]
where $N_{\sa}$ and $N_{\ss}$ are the solid tori from the Dehn surgery. The
trivialization of $w$ over $\partial C$ allows us to extend $w$  to
a line-bundle $w_{\sa} \to Y_{\sa}$ equipped with a trivialization
$\tau_{\sa}$ over
$N_{\sa}$, extending the given trivialization on $\partial C$.  Note
that $w_{\sa}$ is globally trivial on the homology $3$-sphere $Y_{\sa}$,
but the global trivialization differs from $\tau_{\sa}$ on the curve
$b\subset \partial C$ by a map $b\to S^{1}$ of degree $1$. This is
because there is a surface $\Sigma \subset C$ with boundary $b$, and
the original trivialization $\tau$ does not extend over
$\Sigma$. The same remarks apply to $Y_{\ss}$.

On the manifold $Y_{\ss}$, in addition to constructing $w_{\ss}$ as
above, we construct a different line bundle $\tilde w_{\ss} \to
Y_{\ss}$ as follows. Let $\tilde\tau$ be the trivialization of
$w|_{\partial C}$ with the property that $\tilde\tau \tau^{-1}$ is a
map $\partial C\to S^{1}$ with degree $q$ on $b$ and degree $0$ on
$a$. Let $\tilde w_{\ss}$ be obtained by extending $w$ as a trivial
bundle over $N_{\ss}$ extending the trivialization $\tilde\tau$.

If $p$ is odd, then $Y_{\ss}$ has $H^{2}(Y_{\ss};\Z/2) = 0$. When $p$
is even, the construction of $\tilde w_{\ss}$ makes $c_{1}(\tilde
w_{\ss})$ divisible by $2$. So in either case, elements of
$\Rep^{\tilde w_{\ss}}(Y_{\ss})$ correspond to homomorphisms $\rho :
\pi_{1}(Y_{\ss}) \to \SU(2)$.

The following lemma is
straightforward.

\begin{lem}\label{lem:one-to-one-r}
        Restriction to $C$ gives identifications
        \[
        \begin{aligned}
        \Rep^{w_{\sa}} (Y_{\sa}) &\to
                 \Rep^{w}\bigl( \,C \mid \alpha=0
                       \,\bigr) \\
       \Rep^{w_{\ss}} (Y_{\ss}) &\to
                 \Rep^{w}\bigl( \,C \mid p\alpha + q\beta=0
                       \,\bigr)  \\
       \Rep^{\tilde w_{\ss}} (Y_{\ss}) &\to
                 \Rep^{w}\bigl( \,C \mid p\alpha + q\beta= q\pi
                       \,\bigr)                                              
                       \end{aligned}
        \]
        \qed
\end{lem}

The manifold $C$ has $H_{1}(C;\Z) = \Z$, so the representation variety
$\Rep^{w}(C)$ contains reducibles. The next lemma describes their
$\alpha$ and $\beta$ parameters.

\begin{lem}\label{lem:reducible-line}
        If\/ $[A]$ is a reducible element of $\Rep^{w}(C)$, then
        $(\alpha(A), \beta(A))$ lies on the line $\beta = \pi$ mod
        $2\pi\Z$.
\end{lem}

\begin{proof}
        If $[A]$ is a reducible element of $\Rep^{w}(C)$, then $A$ is
        a flat $\SO(3)$ connection on $C$ with cyclic holonomy. The
        holonomy around $b$ is the identity element of $\SO(3)$
        because $b$ bounds the surface $\Sigma$ in $C$. So the
        corresponding $\SU(2)$ connection $\tilde A$ on $E|_{b}$ 
        (regarding $E|_{b}$ as an $\SU(2)$ bundle using $\tau$) has
        holonomy $\pm 1$ in $\SU(2)$. It follows that $\beta$ is $0$
        or $\pi$ mod $2\pi$. We can equip $w$ on $C$ with a connection
        $\theta$ which respects the trivialization $\tau$ on $\partial
        C$ and whose curvature $F_{\theta}$ integrates to $-2\pi i$ on
        $\Sigma$. The $\SU(2)$ connection $\tilde A$ can be uniquely
        extended to a $U(2)$ connection $\tilde A$ on all of $E|_{C}$,
        in such a way that the associated $\SO(3)$ connection is $A$
        and such that the induced connection on $\det(E) =w$ is
        $\theta$. The connection reduces $E$ to a sum of line bundles,
        both of which have curvature $F_{\theta}/2$. The holonomy
        of these line bundles on $b$ is given by
        \[
               \exp \int_{\Sigma}  (F_{\theta}/2) = -1.
        \]
        So $\beta= \pi$ mod $2\pi$ as claimed. This completes the
        proof of the lemma.
\end{proof}

If we suppose that the homology-sphere
$Y_{\sa}$ has a fundamental group with no non-trivial
homomorphisms to $\SU(2)$, then
$\Rep^{w_{\sa}}(Y_{\sa})$ consists of a single reducible element.  By the
previous two lemmas, the $\alpha$ and $\beta$ parameters of this
connection lie on the two line $\alpha = 0$ and $\beta=\pi$. So it is
the point
\[
            v_{\sa} = (0,\pi)
\]
mod $2\pi \Z^{2}$.
Similarly the $\alpha$ and $\beta$ parameters of the reducible
elements in $\Rep^{\tilde w_{\ss}}(Y_{\ss})$ lie on the line $p\alpha + q\beta =
\pi$ mod $2\pi$ and the line $\beta=\pi$. So they are represented by
the points
\[
           v_{\ss,k} = (2\pi k /p,\; \pi)
\]
mod $2\pi\Z^{2}$.
The next lemma is a standard result, from \cite{Floer2} of
\cite{Braam-Donaldson}. We supply the proof for completeness.

\begin{lem}\label{lem:Wstar}
        Suppose $\pi_{1}(Y_{\sa})$ admits no non-trivial homomorphisms
        to $\SU(2)$. For any neighborhood $W$ of $(0,\pi)$, let us write
        \[
                 W^{*} = W \cap \{ \beta \ne \pi \}.
        \]
        Then there exists a symmetric neighborhood $W$ of $(0,\pi)$ such
        that
        \[
              \Rep^{w}( C \mid W^{*} ) = \emptyset.
        \]
\end{lem}

\begin{proof}
        The space $\Rep^{w_{\sa}}(Y_{\sa})$ consists of a single
        point, represented by the $\SO(3)$ connection $A_{\sa}$ with
        trivial holonomy. By the one-to-one correspondence from
        Lemma~\ref{lem:one-to-one-r}, it follows that $\Rep^{w}(C\mid
        (0,\pi))$ consists of a single point $[A]$ represented by an
        $\SO(3)$ connection which trivializes $\g_{E}$.
        We need only show that a neighborhood of $[A]$ in
        $\Rep^{w}(C)$ consists entirely of reducibles. Equivalently,
        writing $\pi$ for $\pi_{1}(C)$, we can study a neighborhood of
        the trivial homomorphism $\rho_{1} : \pi \to \SO(3)$ and show
        that it consists of reducible connections.

        The deformations of $\rho_{1}$ are governed by $H^{1}(\pi;
        \R^{3}) = H^{1}(C) \otimes \R^{3}$, which is a copy of
        $\R^{3}$. It will be sufficient to exhibit a 1-parameter
        deformation of $\rho$ realizing any given vector in this
        $H^{1}$ as its tangent vector and consisting entirely of
        reducibles. This is straightforward. Given $\xi \in \so(3)$, we
        can consider the $1$-parameter family of connections in the
        trivial $\SO(3)$ bundle given by the connection $1$-forms $t
        \xi\eta$, where $\eta$ is a closed $1$-form with period $1$ on
        $C$ and $t\in \R$.
\end{proof}

We need one more lemma before completing the proof of
Proposition~\ref{prop:1/q-triangle}.

\begin{lem}\label{lem:alpha-pi}
       For any $S$, there is a one-to-one correspondence between
       $\Rep^{w}(C\mid S)$ and $\Rep^{w}(C\mid S')$, where $S'$ is the
       translate $S + (\pi,0)$.
\end{lem}

\begin{proof}
        Let $\epsilon$ be an automorphism of the $U(2)$ bundle $E\to
        C$ whose determinant is a function $C \to S^{1}$ which has
        degree $1$ on the curve $a$. (The automorphism $\epsilon$ does
        not belong to the gauge group $\G$, because elements of $\G$
        have determinant $1$.) The element $\epsilon$ acts on the
        space of flat connections $A$ in $\conf(C)$, and gives rise to
        a bijective self-map of the space $\Rep^{w}(C)$:
        \[
             \bar{\epsilon} : \Rep^{w}(C) \to \Rep^{w}(C).
        \]
         This map restricts to a bijection $\bar{\epsilon} :
         \Rep^{w}(C\mid S)\to \Rep^{w}(C\mid S')$. 
\end{proof}

We can now conclude the proof of the proposition. Suppose that
$\pi_{1}(Y_{\sa})$ admits only the trivial
homomorphism to $\SU(2)$, and that the only homomorphisms
$\rho:\pi_{1}(Y_{\ss}) \to \SU(2)$ are those with
cyclic image.
Let $L \subset \R^{2}$ be the
closed line segment
\[
           L = \{ \, (\alpha,\beta) \mid \alpha = 0,  -\pi \le
           \beta\le \pi \,\}
\]
and let $L^{*}$ be the open line-segment obtained by removing the
endpoints. Let $L_{\pi}^{*}$ and $L_{-\pi}^{*}$ be the translates of
this line segment by the vectors $(\pi,0)$ and $(-\pi,0)$.  By
Lemmas~\ref{lem:one-to-one-r} and \ref{lem:reducible-line}, the
hypothesis on $\pi_{1}(Y_{\sa})$ means that
\[
        \Rep^{w}( C \mid L^{*}) = \emptyset.
\]
By Lemma~\ref{lem:alpha-pi},  we therefore have
\[
         \Rep^{w}( C \mid L^{*}_{\pm\pi}) = \emptyset.
\]
Let $P_{1}$ be the line
\[
           P = \{ \, p\alpha+ q\beta= q\pi \,\}
\]
and let $P_{2} = P_{1} - (0,2\pi)$.  The hypothesis on $\pi_{1}(Y_{\ss})$ means
that
$\Rep^{w}( C \mid P_{i})$ consists only of reducibles, lying over the
points on $P_{i}$ where $\beta = \pi$ mod $2\pi$.
Let $S \subset \R^{2}$ be the piecewise-linear arc with vertices at
the points
\[
\begin{aligned}
z_{1} &= (-\pi, 0) \\
z_{2} &= (-\pi, \; -(1-p/q)\pi) \\
z_{3} &= (0,-\pi) \\
z_{4} &= (0,\pi) \\
z_{5} &= (\pi, \; (1-p/q)\pi) \\
z_{6} &= (\pi,0).
          \end{aligned}
\]
Figure~\ref{fig:Dehn} shows the set $S$ in the case $p/q  = 5/3$.
Because $p/q \le 2$, the set is contained
in the region $-\pi \le \beta\le \pi$. If $p/q=2$, then $S$ has four
points on the lines $\beta = \pm\pi$; otherwise it has just two,

Let $S^{*}$ be the complement in $S$ in of the points whose
$\beta$ coordinates are $\pm \pi$.
Given any symmetric neighborhood $U$ of $S$, let $U^{*}$ similarly stand for
\begin{equation}\label{eq:Ustardef}
              U^{*} =  U \setminus \{\, \beta = \pm\pi\, \}.
\end{equation}
We know that
$
         \Rep^{w}(C \mid S^{*}) = \emptyset
$, because $S$ is entirely contained in the union of $L$, $L_{\pm\pi}$
and the two lines $P_{1}$, $P_{2}$.
From Lemma~\ref{lem:Wstar} and the compactness of $\Rep^{w}(C)$, it
follows that there is a symmetric neighborhood $U$ of $S$ such that
\begin{equation}\label{eq:U-empty}
         \Rep^{w}(C \mid U^{*}) = \emptyset.
\end{equation}

\begin{figure}
\begin{center}
\includegraphics[height=2in]{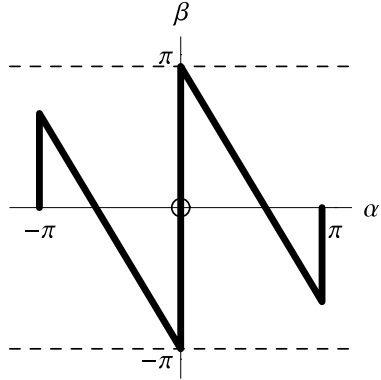}
\end{center}
\caption{\label{fig:Dehn}
The set $S$, for $p/q=5/3$. The $(\alpha,\beta)$ parameters of
reducible elements of $\Rep^{w}(C)$ lie on the dashed lines.}
\end{figure}

We now observe that, given any neighborhood $U$ of $S$, we can find a
smooth
odd function $g$ with period $2\pi$ such that the graph of $-g$ on the
interval $[-\pi,\pi]$ is entirely contained in $U^{*}$.  By
Lemma~\ref{lem:one-to-one-phi} and Lemma~\ref{lem:g}, there exists a
$\phi$ such that
\[
     \Rep^{w}_{\inc,\phi}(Y) = \Rep^{w}\bigl ( \,C \mid \beta =
     -g(\alpha) \,\bigr).
\]
The right hand side is empty because it is contained in the empty set
\eqref{eq:U-empty}. This finishes the proof of the proposition.
\end{proof}

We can reformulate the result of Proposition~\ref{prop:1/q-triangle}
in the special case that $Y_{\sa}$ is $S^{3}$ as follows.

\begin{cor}\label{cor:emptyq}
        Let $K$ be a knot in $S^{3}$ and let $Y_{r}$ be the manifold
        obtained by Dehn surgery with coefficient $r\in \Q$.
        Let $Y_{0}$ be the manifold obtained by
        $0$-surgery, and let $w \to Y_{0}$ be a line bundle whose first
        Chern class is a generator of $H^{2}(Y_{0};\Z)$.
        Suppose $\pi_{1}(Y_{r})$ admits no homomorphism $\rho$ to
        $\SU(2)$ with non-cyclic image.
        Then, if\/ $0 < r < 2$, the manifold $Y_{0}$ admits
        a holonomy deformation $(\inc,\phi)$ so that
        $\Rep^{w}_{\inc,\phi}(Y_{0})$ is empty. \qed
\end{cor}

\section{Proof of the theorem}

\subsection{A stretching argument}

Let $X$ be a closed, oriented $4$-manifold containing a connected,
separating $3$-manifold $Y$. Let $g_{1}$ be metric on $X$ that is
cylindrical on a collar region $[-1,1] \times Y$ containing $Y$ in
$X$. For $L > 0$, let $X_{L}\cong X$ be the manifold obtained from
$X$ by removing the piece $[-1,1]\times Y$ and replacing it with
$[-L,L]\times Y$. There is a metric $g_{L}$ on $X_{L}$ that contains a
cylindrical region of length of $2L$ and agrees
with the original metric on the complement of the cylindrical piece.

Let $v\to X$ be a line bundle, let $E\to X$ be a unitary rank-2 bundle
with $\det(E) = v$, and form the configuration space
$\bonf^{v}(X,E)$ of connections in $\g_{E}$ modulo determinant-1 gauge
transformations of $E$, as we did in the $3$-dimensional case.
In dimension $4$, the bundle $E$ is not determined up to isomorphism
by $v$ alone, so we include it in our notation.
Inside $\bonf^{v}(X,E)$ is the
moduli space of anti-self-dual connections,
\[
       M^{v}(X,E) = \{ \, [A] \in \bonf^{v}(X,E) \mid  F^{+}_{A} = 0
       \, \}.
\]
For each $L > 0$, we also have a moduli space \[ M^{v}(X_{L},E)
\subset \bonf^{v}(X_{L},E). \] (We do
not take the trouble to introduce the additional notation $v_{L}$ and
$E_{L}$ for the corresponding bundles on $X_{L}$.)

Let $(\inc,\phi)$ be data for a holonomy perturbation for the bundle
$E|_{Y}$. Following \cite{FloerOriginal, Floer2, Donaldson-Book}, we
shall use $\phi$ also to perturb the anti-self-duality equations on $X_{L}$.
We use $\inc$ to embed $[-L,L]\times S^{1} \times D$ into $X_{L}$, and
let $\mu_{X}$ be the $2$-form on the cylindrical part $[-L,L] \times Y$
obtained by pulling back $\mu$ from $D$ and pushing forward using this
embedding. We choose a trivialization of $v = \det(E)$ on the image of
the embedding so that each $\SO(3)$ connection in $\g_{E}$ determines
uniquely an $\SU(2)$ connection. For each $A$, the
holonomy around the circles
defines, as before, a section $H_{A}$ over $[-L,L] \times \im(\inc)$
of the bundle $\mathrm{Aut}(E)$, and we obtain
\[
               \phi'(H_{A}) \in  C^\infty([-L,L]\times\im(\inc);
               \g_{E}).
\]
For $L>1$, let $\beta : X_{L} \to [0,1]$ be a smooth cut-off function,
supported in $[-L,L]\times Y$ and equal to $1$ on $[-L+1, L-1]\times
Y$.  On $X_{L}$, the \emph{perturbed anti-self-duality equation} is
the equation
\begin{equation}
    \label{eq:pertASD}
                F^{+}_{A} + \beta \phi'(H_{A}) \mu^{+} = 0.
\end{equation}
We define the corresponding moduli space:
\begin{equation}
    M^{v}_{\phi}(X_{L},E) = \{ \, [A] \in \bonf^{v}(X_{L},E) \mid
    \text{equation \eqref{eq:pertASD} holds} \,\}.
\end{equation}

\begin{prop}
    Let $w = v|_{Y}$. Suppose that there is a holonomy perturbation
    on $Y$ such that the perturbed representation variety
    $\Rep^{w}_{\inc,\phi}(Y)$ is empty. Then for each $E$ with
    determinant $v$ on $X$, there exists an $L_{0}$
    such that $M^{v}_{\phi}(X_{L},E)$ is also empty,
    for all $L\ge L_{0}$.
\end{prop}

\begin{proof}
    The proof is some subset of a standard discussion of holonomy
    perturbations and compactness in Floer homology theory (see
    \cite{Floer2,Braam-Donaldson,Donaldson-Book}).
    Suppose on the contrary that we can find $[A_{i}]$  in
    $M^{v}_{\phi}(X_{L_{i}},E)$ for an increasing, unbounded sequence
    of lengths $L_{i}$. We start as usual with the fact that the
    quantity
    \[
    \begin{aligned}
    \cE(A_{i}) &= \int_{X_{L_{i}}} \tr (F_{A_{i}} \wedge F_{A_{i}}) \\
             &= \| F^{-}_{A_{i}} \|^{2} - \| F^{+}_{A_{i}} \|^{2}
             \end{aligned}
    \]
    is independent of $i$ and depends only on the Chern numbers of the
    bundle $E$. (The norms are $L^{2}$ norms.)
    We write this quantity as the sum of three terms:
    \[
           \cE(A_{i}) = \cE(A_{i}\mid X^{1}) + \cE(A_{i}\mid X^{2})
                         + \cE(A_{i}\mid X^{3}_{i}), 
    \]
    where
    \[
  \begin{aligned}
    X^{1}  &= X_{L_{i}} \setminus \bigl([-L_{i},L_{i}] \times Y \bigr)\\
    X^{2}  &= \bigl( [-L_{i},-L_{i}+1] \times Y \bigr) \cup  \bigl(
    [L_{i}-1,L_{i}] \times Y \bigr)  \\
    X^{3}_{i}  &= [-L_{i}+1,L_{i}-1] \times Y .
  \end{aligned}
    \]
Only the third piece has a geometry which depends on $i$.  From the
equation \eqref{eq:pertASD}, we have
\[
           \cE(A_{i}\mid X^{1})  \ge 0
\]
because $\beta$ is zero on $X^{1}$. The second term in equation
\eqref{eq:pertASD} is pointwise uniformly bounded, so
\[
           \cE(A_{i}\mid X^{2})  \ge - C_{2}
\]
where $C_{2}$ is independent of $i$. Because the sum of the three terms is
constant, we deduce that
\[
            \cE(A_{i}\mid X^{3}_{i})  \le K,
\]
where $K$ is independent of $i$.

To understand the term $\cE(A_{i}\mid X^{3}_{i})$ better, one must
reinterpret \eqref{eq:pertASD}.  On $X^{3}_{i}$, the function $\beta$
is $1$. Identify $E$ on this cylinder with the pull-back of a bundle
$E_{Y} \to Y$, and choose a gauge representative $A_{i}$ for $[A_{i}]$ in
temporal gauge. Write
\[
             A_{i}(t) = A_{i}|_{\{t\}\times Y}, \quad (-L_{i}+1 \le t \le
             L_{i}-1).
\]
Thus $A_{i}(t)$ becomes a path in the space of connections $\conf(Y;E_{Y})$.  The
equation \eqref{eq:pertASD} is equivalent on $X^{3}_{i}$ to the
condition that $A_{i}(t)$ solves the downward gradient flow equation
for the perturbed Chern-Simons functional on $\conf(Y;E_{Y})$:
\[
             \frac{d}{dt} A_{i}(t) = - \grad ( \CS + \Phi ).
\]
In particular, $\CS + \Phi$ is monotone decreasing along the path (or
constant). The function $|\Phi|$ is a bounded function on
$\conf(Y;E_{Y})$: we can write
\[
        |\Phi| \le K'.
\]
The change in $\CS$ is equal to the quantity $-\cE$:  that is,
\[
\begin{aligned}
\CS\bigr(A_{i}(-L_{i} + 1)\bigl) -
            \CS\bigr(A_{i}(L_{i} - 1)\bigl)  &=  \cE(A_{i}\mid
            X^{3}_{i})  \\
            &\le K
            \end{aligned}
\]
So from the bound on $|\Phi|$ we obtain
\[
(\CS+\Phi)\bigr(A_{i}(-L_{i} + 1)\bigl) -
            (\CS+\Phi)\bigr(A_{i}(L_{i} - 1)\bigl)   \le K +  2K'. 
\]

Now let $\delta>0$ be given.  Because $\CS + \Phi$ is decreasing and
the total drop is bounded by $K + 2K'$, we can find intervals
\[
       (a_{i}, b_{i}) \subset [-L_{i}+1, L_{i}+1]
\]
of length $\delta$, so that the drop in $\CS + \Phi$ along $(a_{i},
b_{i})$ tends to zero as $i$ goes to infinity. Because the equation is
a gradient-flow equation, this means
\[
           \lim_{i\to \infty} \int_{a_{i}}^{b_{i}} \| \grad ( \CS +
           \Phi )(A_{i}(t)) \|_{L^{2}(Y)}^{2}\,dt = 0.
\]
We have an expression for $\grad\Phi$ as a uniformly bounded form, so
\[
              \limsup_{i\to\infty}
              \int_{a_{i}}^{b_{i}} \|F_{A_{i}(t)} \|_{L^{2}(Y)}^{2}\,dt
              \le\delta J
\]
for some constant $J$ depending on $\phi$.
So given any $\epsilon > 0$, we can find a $\delta>0$ and a sequence
of intervals $(a_{i},b_{i})$  of length $\delta$ so that
\[
         \int_{(a_{i},b_{i})\times Y}  | F_{A_{i}} |^{2}
         \,d\mathrm{vol} \le \epsilon
\]
for all $i\ge i_{0}$.
We now regard the $A_{i}$ as connections on the fixed cylinder
$(0,\delta)\times Y$.
At this point, if $\epsilon$ is smaller than the
threshold for Uhlenbeck's gauge fixing theorem on the $4$-ball, we can
find $4$-dimensional gauge transformations on the cylinder
so that, after applying these gauge transformations
and passing to a subsequence, the connections converge in $C^{\infty}$
on compact subsets. (See for example \cite[section 5.5]{Donaldson-Book}.)

If $A$ is the limiting connection on $(0,\delta)\times Y$, in temporal
gauge, then the
function $\CS + \Phi$ is constant along the path $A(t)$. It follows
that $A(t)$ is constant and is
a critical point of $\CS+\Phi$. This tells us that
$[A(t)]$ belongs to the perturbed representation variety
$\Rep^{w}_{\inc,\phi}(Y)$, which we were supposing to be empty.
\end{proof}

The proposition above has the following corollary for the Donaldson
polynomial invariants. (Our notation and conventions for these
invariants is taken from \cite{KM-structure}.)

\begin{cor}\label{cor:Dvanishes}
    Let $X$ be an admissible $4$-manifold in the sense of
    \cite{KM-structure}, so that its Donaldson polynomial invariants
    $D_{X}^{v}$ are defined. (For example, suppose $H_{1}(X;\Z)$ is
    zero and $b^{+}(X)$ is greater than $1$.) Then, under the
    assumptions of the previous proposition, the polynomial invariants
    are identically zero, regarded as a map
    \[
               D_{X}^{v} : \A(X) \to \Z.
    \] 
\end{cor}

\begin{proof}
    The definition of $D_{X}^{v}$ involves first choosing a Riemannian
    metric on $X$ so that the moduli spaces $M^{v}(X,E)$ are smooth
    submanifolds of $\bonf^{v}(X,E)$, containing no reducibles and cut
    out transversely by the equations. If $X$ is admissible, then this
    can always be done, by changing the metric inside a ball in $X$.
    The value of the invariant is then defined as a signed count of
    the intersection points between $M^{v}(X,E)$ and some
    specially-constructed finite-codimension submanifolds of
    $\bonf^{v}(X,E)$.  This part of the construction of $D^{v}_{X}$
    involves only transversality arguments, which can be carried out
    equally with $M^{v}_{\phi}(X_{L},E)$ in place of $M^{v}(X,E)$, for
    any fixed $L$.  That the signed count is independent of the
    choices made, in the unperturbed setting, is a consequence of the
    compactness theorem for the moduli space. The Uhlenbeck
    compactification works the same way for $M^{v}_{\phi}(X_{L},E)$ as it
    does for the unperturbed anti-self-duality equations (see
    \cite{Donaldson-Book} for example); so the Donaldson invariants
    can be defined using the perturbed moduli spaces. Each moduli
    space is empty once $L$ is large enough, so the invariants are
    zero.
\end{proof}

\subsection{Concluding the proof}

The rest of the argument is essentially the same as the proof of the
main theorem in \cite{KM-PP}. 
Let $K$ be a knot in $S^{3}$ that is a counterexample to
Theorem~\ref{thm:hRP3}.  We will obtain a
contradiction.

The manifold $Y_{0}$ obtained by zero-surgery admits a taut foliation
and is not $S^{1}\times S^{2}$, by the results of \cite{Gabai3}. The
following proposition is proved in \cite{KM-PP} using the results of
\cite{Eliashberg-Thurston} and \cite{Eliashberg}:

\begin{prop}\label{prop:Summary}
    Let $Y$ be a closed orientable $3$-manifold admitting an oriented taut
    foliation. Suppose $Y$ is not $S^{1}\times S^{2}$. Then $Y$ can be
    embedded as a separating hypersurface in a closed symplectic $4$-manifold
    $(X,\Omega)$. Moreover, we can arrange that $X$ satisfies the
    following additional conditions.
    \begin{enumerate}
        \item The first homology $H_{1}(X;\Z)$ vanishes.

        \item The euler number and signature of $X$ are the same as
        those of some smooth hypersurface in
        $\CP^{3}$, whose degree is even and not less than $6$.

         \item The restriction map $H^{2}(X;\Z) \to H^{2}(Y;\Z)$ is
        surjective. \label{it:H2Onto}

        \item The manifold $X$ contains a tight surface of positive
        self-intersection number, and a sphere of self-intersection
        $-1$. \qed
    \end{enumerate}
\end{prop}

We apply this proposition to the manifold $Y_{0}$, to obtain an $X$
with all of the above properties.
Using the results of \cite{Feehan-Leness}, it was shown in
\cite{KM-PP} that a $4$-manifold satisfying these conditions satisfies
Witten's conjecture relating the Seiberg-Witten and Donaldson
invariants. (See \cite[Conjecture 5 and Corollary 7]{KM-PP} for an
appropriate statement of Witten's conjecture in this context.) Because
$X$ is symplectic, its Seiberg-Witten invariants are non-trivial by
\cite{Taubes1}. For the same reason, $X$ has Seiberg-Witten simple
type. From Witten's conjecture, it follows that the Donaldson
invariants $D^{v}_{X}$ are non-trivial, for all $v$ on $X$.

By the penultimate condition on $X$ in Proposition~\ref{prop:Summary},
we can choose $v\to X$ so that $c_{1}(v)$ restricts to a generator of
$H^{2}(Y_{0};\Z)$. Write $w = v|_{Y_{0}}$.  If $K$ is a counterexample
to Theorem~\ref{thm:hRP3}, then Corollary~\ref{cor:emptyq} tells us
there is a holonomy perturbation $\phi$ such that
\[
       \Rep^{w}_{\inc,\phi}(Y_{0}) = \emptyset.
\]
Corollary~\ref{cor:Dvanishes} then tells us that
$D^{v}_{X}$ is zero. This is the contradiction. \qed

\subsection{Further remarks}

An analysis of the proof of
Theorem~\ref{thm:hRP3} reveals that it
proves a slightly stronger result (stronger, that is, if one is
granted the results of \cite{Gabai3}). For example, we can state:

\begin{thm}
    Let $N$ be an embedded solid torus in an irreducible closed
    $3$-manifold $Y$ with $H_{1}(Y) = \Z$. Let $C = Y \setminus
    N^{\circ}$ be the complementary manifold with torus boundary.

    Then there is at most one Dehn filling of $C$ which yields a
    homotopy sphere. Indeed, for all but one slope, the fundamental
    group of the manifold obtained by Dehn filling admits a
    non-trivial homomorphism to $\SU(2)$. \qed
\end{thm}

The point here is that the original hypothesis need not be that $K$ is
a non-trivial knot in $S^{3}$. What one wants is that zero-surgery on
$K$ should be an irreducible homology $S^{1}\times S^{2}$; and if we
make this our hypothesis, then we can also consider the case that $K$
is a knot in (for example) a homotopy sphere.

One can also ask whether there is a non-trivial extension of
Theorem~\ref{thm:hRP3} to other integer surgeries.
The results of \cite{KMOS} show that surgery with coefficient $3$ or
$4$ on a non-trivial knot cannot be a lens space. It
would be interesting to know whether the fundamental groups of $Y_{3}$
and $Y_{4}$ must admit homomorphisms to $\SU(2)$ with non-abelian
image when $K$ is non-trivial. Surgery with coefficient $+5$ on the
right-handed trefoil produces a lens space, so
one does not expect to extend Theorem~\ref{thm:hRP3} further in the
direction of integer surgeries without additional hypotheses.
Dunfield \cite{Dunfield} has provided an example of a non-trivial knot in
$S^{3}$ for
which the Dehn filling $Y_{37/2}$ has a fundamental group which is not
cyclic but admits no homomorphism to $\SU(2)$ (or even $\SO(3)$)
with non-abelian image.
(The knot is the $(-2,3,7)$ pretzel knot,
for which 
$Y_{18}$ and $Y_{19}$ are both lens spaces \cite{Fintushel-Stern}.) This
example shows that the property of having cyclic fundamental group and the
property of admiting no cyclic homomorphic image in $\SU(2)$ are in
general different for $3$-manifolds obtained by Dehn surgery.

\bibliography{hrpthree}

\end{document}